\documentclass[english]{article}
\usepackage{preamble}

\begin{document}

\title{A Remark on the Number of Maximal Abelian Subgroups}
\maketitle

\begin{abstract}
    The number of maximal abelian subgroups of a finite $p$-group is shown to be congruent to $1$ modulo $p$. 
\end{abstract}

We say that a subgroup of a group $G$ is \textit{maximal abelian}, if it is abelian and not  properly contained in any larger abelian subgroup of $G$. That is, if it is maximal with respect to inclusion among abelian subgroups of $G$. In particular, the maximal abelian subgroups of $G$ need not be maximal subgroups of $G$ and may have different orders. The purpose of this note is to show the following counting result:  

\begin{thm}\label{thm}
    Let $G$ be a finite $p$-group. The number of maximal abelian subgroups of $G$ is congruent to $1$ modulo $p$.
\end{thm}

While there exist numerous counting results of similar flavour (many of which can be found in \cite{berkovich2011groups}), this particular one seems to have not been previously observed. The proof employs the standard technique of M{\"o}bius inversion on the subgroup lattice.  
 
\begin{proof}
    To facilitate induction, we shall prove, more generally, that for every abelian subgroup $H\le G$, the number 
    \[
        g_G(H) = |\{A\le G \mid H\sseq A \text{ and }  A \text{ is maximal abelian}\}|
    \] 
    is congruent to $1$ modulo $p$. The theorem follows by considering the trivial subgroup $H=\{1\}\le G$. We shall prove the claim by induction on $[G:H]$. The base of the induction is $[G:H]=1$, where we have $H=G$ and hence $g_G(H)=1$.
    
    For the inductive step, we begin by reducing to the case $H=Z$, where $Z$ is the \textit{center} of $G$.
    First, note that every abelian subgroup of $G$ that contains $H$ must lie in the centralizer subgroup $C(H)$ of $H$ in $G$. Hence, 
    \(
        g_G(H) = g_{C(H)}(H).
    \)
    Therefore, if $C(H)\subsetneq G$, then we are done by the inductive hypothesis. It thus suffices to consider only the case $C(H)=G$, or equivalently,  $H\sseq Z$. Second, note that every maximal abelian subgroup of $G$ must contain $Z$, so we get 
    \(
        g_G(H) = g_G(Z).
    \) 
    Therefore, if $H\subsetneq Z$, then we are once again done by the inductive hypothesis. It remains to consider the case $H = Z$.
    
    Let $\cS(G)$ be the inclusion lattice of (all) subgroups of $G$. The M{\"o}bius function for this lattice is given by (see \cite[Theorem 2]{weisner1935some} or \cite{hall1936eulerian})
    \[
        \mu(S,T) = 
            \begin{cases}
                (-1)^k p^{\binom{k}{2}} & \text{if } S\triangleleft T \text{ and } T/S\simeq (\ZZ/p)^k \\
                0 & \text{ else.}
            \end{cases}
    \]
    Now, let $f_G \colon \cS(G) \to \ZZ$ be the indicator function of the subset of maximal abelian subgroup of $G$. Observe that $g_G$ is the accumulative function of $f_G$ in the sense that
    \[
        g_G(S) = \sum_{S\sseq T} f_G(T).
    \]
    We want to prove that $g_G(Z)$ is congruent to $1$ modulo $p$. If $Z$ is itself a maximal abelian subgroup of $G$, then $g_G(Z)=1$ and we are done. Otherwise, by M{\"o}bius inversion, we have
    \[
        f_G(Z) = \sum_{Z\sseq T} \mu(Z,T)g_G(T).
    \]
    Since, by assumption, $Z$ is not maximal abelian, we get $f_G(Z)=0$. Isolating the term $T=Z$ in the sum, we obtain 
    \[
        0 =  g_G(Z) + \sum_{Z\subsetneq T} \mu(Z,T)g_G(T)
    \]
    \[
        g_G(Z) = \sum_{Z\subsetneq T} -\mu(Z,T)g_G(T).
    \]
    We now analyze the terms in the sum modulo $p$. To begin with, if $T$ is not abelian, then $g_G(T)=0$, and if $T$ is abelian, then by the inductive hypothesis $g_G(T)$ is $1$ modulo $p$. thus,
    \[
        g(Z) \equiv \sum_{Z\subsetneq T \text{ abelian}} -\mu(Z,T)
        \quad (\mathrm{mod}\;\,p).
    \]
    Furthermore, by the explicit formula of $\mu$ above, we have that $-\mu(Z,T)$ is zero modulo $p$, unless $T/Z\simeq \ZZ/p$, where $-\mu(Z,T)=1$ (note that $Z$ is always normal in $T$). Hence, $g_G(Z)$ is congruent modulo $p$ to the number of abelian subgroups $T$ of $G$, such that $Z\sseq T$ and $T/Z\simeq \ZZ/p$. It is a standard fact that if a quotient of a group by its center is cyclic, then the group is abelian. Thus, $g(Z)$ is congruent modulo $p$ to the number of order $p$ subgroups of $G/Z$. The claim now follows from the fact that the number of order $p$ subgroups of a non-trivial finite $p$-group is $1$ modulo $p$.
\end{proof}  

\begin{rem}
    As suggested to me by Peter M{\"u}ller, one can avoid the M{\"o}bius inversion formula and proceed instead more directly by a double counting argument. The idea is to count the number of pairs $(c,A)$, where $A\le G$ is a maximal abelian subgroup and $c \in G/Z$ is a coset, such that $c \sseq A$. The details of how to recast the inductive argument above in this perspective are left to the reader. 
\end{rem}
 
\textbf{Acknowledgements. } I wish to thank the Group Pub Forum community for a useful discussion regarding the subject of this note. In particular, I would like to thank Peter M{\"u}ller for useful comments on the proof and David Craven and Avinoam Mann for valuable bibliographical information.

 

\bibliographystyle{alpha}
\phantomsection\addcontentsline{toc}{section}{\refname}
\bibliography{MaxAb}

\end{document}